\documentclass[11pt]{article}  

\usepackage{color}              
\usepackage{graphicx,setspace}
\usepackage[]{amsmath}
\usepackage{amssymb}
\usepackage{fullpage}
\usepackage{enumitem}
\usepackage{hyperref}
\usepackage{amsfonts, amsthm, bbm}
\usepackage{fullpage}

\definecolor{MyDarkBlue}{rgb}{0,0.08,0.50}
\definecolor{BrickRed}{rgb}{0.65,0.08,0}

\hypersetup{
colorlinks=true,       
    linkcolor=MyDarkBlue,          
    citecolor=BrickRed,        
    filecolor=red,      
    urlcolor=cyan           
}

\renewcommand{\geq}{\geqslant}
\renewcommand{\leq}{\leqslant}
\newtheorem{thm}{Theorem}

\newtheorem{cor}[thm]{Corollary}
\newtheorem{prop}[thm]{Proposition}
\newtheorem{lem}[thm]{Lemma}
\theoremstyle{definition}
\newtheorem{rem}[thm]{Remark}

\newcommand\Sc{\mathcal{S}}
\newcommand\Zb{\mathbb{Z}}

\begin{document}

\title{{\Large {\bf The compensation approach for walks with small steps\\in the quarter plane}}}

\author{Ivo J.B.F. Adan\footnotemark[1] \footnotemark[3]  \and
Johan S.H. van Leeuwaarden\footnotemark[1] \footnotemark[3]\and
        Kilian Raschel\footnotemark[2] \footnotemark[3]
        }

\date{\today}

\maketitle

\footnotetext[1]{Department of Mathematics and Computer Science, Eindhoven University of Technology,
P.O.\ Box 513, 5600 MB  Eindhoven, The Netherlands
}

\footnotetext[2]{CNRS and Universit\'e de Tours,
Facult\'e des Sciences et Techniques,
Parc de Grandmont,
37200 Tours, France}

\footnotetext[3]{E-mails:  \url{i.j.b.f.adan@tue.nl}, \url{j.s.h.v.leeuwaarden@tue.nl},
        \url{kilian.raschel@lmpt.univ-tours.fr}}

\begin{abstract}
This paper is the first application of the compensation approach (a well-established theory in probability theory) to counting problems. We discuss how this method can be applied to a general class of walks in the quarter plane $\Zb_{+}^{2}$ with a step set that is a subset of $\{(-1,1),(-1,0),(-1,-1),(0,-1),(1,-1)\}$ in the interior of $\Zb_{+}^{2}$. We derive an explicit expression for the generating function which turns out to be nonholonomic, and which can be used to obtain exact and asymptotic expressions for the counting numbers.
\end{abstract}


\vspace{1.5mm}

\noindent {\it AMS $2000$ Subject Classification}: 05A15, 05A16

\section{Introduction}
\label{sec:intro}
In the field of enumerative combinatorics, counting walks on a lattice is among the most classical topics.
While counting problems have been largely resolved for unrestricted walks on $\Zb^2$, walks that are confined to the quarter plane $\Zb_{+}^{2}$ still pose considerable challenges. In recent years, much progress has been made, in particular for walks in the quarter plane with small steps, which means that the step set $\Sc$ is a subset of $\{(i,j) : |i|,|j|\leq 1\}\setminus \{(0,0)\}$. Bousquet-M\'{e}lou and Mishna \cite{BMM} constructed a thorough classification of these models. By definition, there are $2^{8}$ such models, but after eliminating trivial cases and exploiting equivalences, it is shown in \cite{BMM} that there are $79$ inherently different problems that need to be studied. Let $q_{i,j,k}$ denote the number of paths in $\Zb_{+}^{2}$ having length $k$, starting from $(0,0)$ and ending in $(i,j)$, and define the generating function (GF) as
     \begin{equation}
     \label{CGF}
          Q(x,y;z)=\sum_{i,j,k=0}^{\infty}q_{i,j,k}x^{i}y^{j}z^{k}.
     \end{equation}
There are then two key challenges:
\begin{enumerate}
     \item \label{Challenge_1} Finding an explicit expression for $Q(x,y;z)$.
\item \label{Challenge_2} Determining the nature of $Q(x,y;z)$: is it holonomic
       (the vector space
        over $\mathbb{C}(x,y,z)$---the field of rational functions in the three variables $x,y,z$---spanned by the set of all derivatives of $Q(x,y;z)$ is finite dimensional, see \cite[Appendix B.4]{FLAJ})? And in that event, is it algebraic, or even rational?
\end{enumerate}
The common approach to address these challenges is to start from a functional equation for the GF, which for the walks with small steps takes the form (see~\cite[Section 4]{BMM})
     \begin{equation}
     \label{FE}
          K(x,y;z)Q(x,y;z) =  A(x)Q(x,0;z)+B(y)Q(0,y;z) -\delta Q(0,0;z)-x y/z,
     \end{equation}
where $K$, $A$ and $B$ are polynomials of degree two in $x$ and/or $y$, and $\delta$ is a constant. For $z=1/|\Sc|$, (\ref{FE}) belongs to the generic class of functional equations (arising in the probabilistic context of {\em random} walks) studied and solved in the book \cite{FIM}. For general values of $z$, the analysis of (\ref{FE}) for the $79$ above-mentioned models has been carried out in \cite{RaGe,Ra}, which settled Challenge \ref{Challenge_1}.

In order to describe the results regarding Challenge \ref{Challenge_2}, it is worth to define the \emph{group of the walk}, a notion
introduced by Malyshev \cite{MAL}. This is the group of birational transformations $W=\langle \xi,\eta\rangle$, with
\begin{equation}
          \xi(x,y)= \Bigg(x,\frac{1}{y}\frac{\sum_{(i,-1)\in\Sc}x^{i}}
          {\sum_{(i,+1)\in\mathcal{S}}x^{i}}\Bigg),
          \quad \eta(x,y)=\Bigg(\frac{1}{x}\frac{\sum_{(-1,j)\in\mathcal{S}}y^{j}}
          {\sum_{(+1,j)\in\mathcal{S}}y^{j}},y\Bigg),
     \end{equation}
which leaves invariant the  function $\sum_{(i,j)\in\Sc}x^{i}y^{j}$. Clearly, $\xi\circ\xi=\eta\circ\eta=\text{id}$, and $W$ is a dihedral group of even order larger than or equal to four.

Challenge \ref{Challenge_2} is now resolved for all $79$ problems. It was first solved for the $23$ models that have a finite group. The nature of the GF (as a function of the three variables $x,y,z$) was determined in \cite{BMM} for $22$ of these $23$ models: $19$ models turn out to have a GF that is holonomic but nonalgebraic, while $3$ models have a GF that is algebraic. For the $23$rd model, defined by $\mathcal{S}=\{(-1,0),(-1,-1),(1,0),(1,1)\}$ and known as \emph{Gessel's walk}, it was proven in \cite{BK} that $Q(x,y;z)$ is algebraic. Alternative proofs for the nature of the GF for these $23$ problems were given in \cite{FR}.
For the remaining $56$ models, which all have an infinite group, it was first shown that $5$ of them (called singular) have a nonholonomic GF, see \cite{MMM,MM2}. Bousquet-M\'elou and Mishna \cite{BMM} have conjectured that the $51$ other models also have a nonholonomic GF. Partial evidence was provided in \cite{BK1,Ra}, and the proof was recently given in \cite{KRH}.

\begin{figure}[!ht]
\begin{center}
\begin{picture}(73.00,70.00)
\includegraphics{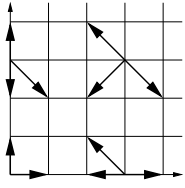}
\end{picture}
\end{center}
\vspace{-3mm}
\caption{The jumps of the walk}
\label{Jumps}
\end{figure}

In this paper we consider walks on $\Zb_{+}^{2}$ with small steps that do not fall into the class considered in \cite{BMM}. The classification in \cite{BMM} builds on the assumption that the steps on the boundaries are the same steps (if possible) as those in the interior. However, when the behavior on the boundaries is allowed to be different, we have many more models to consider. There are several motivations to permit arbitrary jumps on the boundary. First, this allows us to extend the existing results, and thus to have a better understanding of the combinatorics of the walks confined to a quarter plane. Second, some studies in probability theory (see e.g.~\cite{KM}) suggest that the behavior of the walk on the boundary has a deep and interesting influence on the quantities of interest. In this paper we first consider the walk in Figure \ref{Jumps}, with steps taken from $\mathcal{S}=\{(-1,1),(-1,-1),(1,-1)\}$ in the interior, $\mathcal{S}_H=\{(-1,1),(-1,0),(1,0)\}$ on the horizontal boundary, $\mathcal{S}_V=\{(0,1),(0,-1),(1,-1)\}$ on the vertical boundary, and $\mathcal{S}_0=\{(0,1),(1,0)\}$ in state $(0,0)$.

Notice that the step set $\mathcal{S}=\{(-1,1),(-1,-1),(1,-1)\}$ in the framework of \cite{BMM} would render  a trivial walk, since
the walk could never depart from state $(0,0)$. However, by choosing the steps on the boundaries as in $\mathcal{S}_H$, $\mathcal{S}_V$ and $\mathcal{S}_0$, it becomes possible to start walking from state $(0,0)$, and we have a rather intricate counting problem on our hands.

\begin{figure}[!ht]
\begin{center}
\begin{picture}(73.00,70.00)
\includegraphics{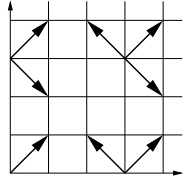}
\end{picture}
\end{center}
\vspace{-3mm}
\caption{Walk considered in \cite{MM2}}
\label{of}
\end{figure}

It turns out that our walk has an infinite group. To see this, observe that the interior step set in Figure \ref{Jumps} and the one represented in Figure \ref{of} have isomorphic groups (\cite[Lemma 2]{BMM} says that two step sets differing by one of the eight symmetries of the square have isomorphic groups), and notice that the group associated with the step set in Figure \ref{of} is infinite (\cite[Section 3.1]{BMM}).

Since our walk has an infinite group, the approach in \cite{BMM} cannot be applied. The methods developed in \cite{MM2} also fail to work. Indeed, the main tool used there is an expression of $Q(x,0;z)$ and $Q(0,y;z)$ as series involving the iterates of the roots of the kernel. While these series are convergent when the transitions $(-1,0),(-1,-1),(0,-1)$ are absent, they become strongly divergent in our case, see \cite[Chapter $6$]{FIM}. The approach via boundary value problems of \cite{FIM} seems to apply, but this is more cumbersome than in \cite{RaGe,Ra}. Indeed,
contrary to the $79$ models studied there, for which in (\ref{FE}), $A$ and $B$ depend on one variable, and $\delta$ is constant, these quantities
are now polynomials in two variables, since the functional equation \eqref{FE} for our walk becomes (for a proof, see Equations \eqref{eq_i>0_j>0_sum}--\eqref{eq_i=0_j=0_sum} in Section \ref{The_compensation_approach})
     \begin{equation}\label{func2}
          K(x,y;z)Q(x,y;z)=
          [1+x^2-x^2 y-y]Q(x,0;z)+[1+y^2-xy^2-x]Q(0,y;z) + [x+y-1]Q(0,0;z)-xy/z
     \end{equation}
with the kernel
\begin{equation}
\label{kern2}
K(x,y;z)=1+x^2+y^2-xy/z.
\end{equation}

We shall derive an explicit expression for $Q(x,y;z)$, which turns out to be a meromorphic function of $x,y$ with infinitely many poles in $x$ and in $y$. This implies that $Q(x,y;z)$ is nonholonomic as a function of $x$, and as a function of $y$. We note that for one variable, a function is holonomic if and only if it satisfies a linear differential equation with polynomial coefficients (see \cite[Appendix B.4]{FLAJ}), so that a holonomic function must have finitely many poles, since the latter are found among the zeros of the polynomial coefficients of the underlying differential equation. Accordingly, the trivariate function $ Q(x,y;z)$ is nonholonomic as well, since the holonomy is stable by specialization of a variable (see \cite{FLAJ}). As  in \cite{KRH,Ra}, this paper therefore illustrates the intimate relation between the infinite group case and nonholonomy.

The technique we are using is the so-called {\it compensation approach}. This technique has been developed in a series of papers \cite{MR1138205,MR1080417,MR1241929} in the probabilistic context of {\it random} walks; see \cite{MR1833660} for an overview. It does not aim directly at obtaining a solution for the GF, but rather tries to find a solution for its coefficients
\begin{equation}
     \label{notation_z}
          q_{i,j}(z)=\sum_{k=0}^{\infty} q_{i,j,k}z^{k}.
\end{equation}
These coefficients satisfy certain recursion relations, which differ depending on whether the state $(i,j)$ lies on the boundary or not.
The idea is then to express $q_{i,j}(z)$ as a linear combination of products $\alpha^i\beta^j$, for pairs $(\alpha,\beta)$ such that
\begin{equation}  \label{dasas}
          K(1/\alpha,1/\beta;z)=0.
\end{equation}
By choosing only pairs $(\alpha,\beta)$ for which \eqref{dasas} is satisfied, the recursion relations for $q_{i,j}(z)$ in the interior of the quarter plane are satisfied by  any linear combination of products $\alpha^i\beta^j$, by virtue of the linearity of the recursion relations. The products have to be chosen such that the recursion relations on the boundaries are satisfied as well. As it turns out, this can be done by alternatingly compensating for the errors on the two boundaries, which eventually leads to an infinite series of products.

This paper is organized as follows. In Section \ref{The_compensation_approach} we obtain an explicit expression for the generating function $Q(x,y;z)$ by applying the compensation approach.
In Section
\ref{Asymptotic_analysis} we derive an asymptotic expression for the coefficients $q_{i,j,k}$ for large values of $k$,
using the technique of singularity analysis. Because this paper is the first application of the compensation approach to counting problems, we also shortly discuss in Section \ref{Discussion} for which classes of walks this compensation approach might work.

\section{The compensation approach}
\label{The_compensation_approach}
We start from the following recursion relations:
     \begin{align}
          q_{i,j,k+1}&=q_{i-1,j+1,k}+q_{i+1,j-1,k}+q_{i+1,j+1,k},  &i,j&\geq 1,k\geq 0, \label{eq_i>0_j>0}\\
          q_{i,0,k+1}&=q_{i-1,1,k}+q_{i+1,1,k}+q_{i-1,0,k}+q_{i+1,0,k},  &i&\geq 1,k\geq 0,\label{eq_i>0_j=0}\\
          q_{0,j,k+1}&=q_{1,j-1,k}+q_{1,j+1,k}+q_{0,j-1,k}+q_{0,j+1,k}, &j&\geq1 ,k\geq 0,\label{eq_i=0_j>0}\\
          q_{0,0,k+1}&=q_{0,1,k}+q_{1,1,k}+q_{1,0,k}, & & k\geq 0.\label{eq_i=0_j=0}
     \end{align}
Since $q_{i,j,0}=0$ if $i+j>0$ and $q_{0,0,0}=1$, these relations uniquely determine all the counting numbers $q_{i,j,k}$.
Multiplying the relations (\ref{eq_i>0_j>0})--(\ref{eq_i=0_j=0}) by $z^k$ and summing w.r.t.\ $k\geq 0$
leads to (with the generating functions $q_{i,j}= q_{i,j}(z)$ defined in \eqref{notation_z})
     \begin{align}
          \label{eq_i>0_j>0_sum}
          q_{i,j}/z&=q_{i-1,j+1}+q_{i+1,j-1}+q_{i+1,j+1},&i,j&\geq 1,\\
          \label{eq_i>0_j=0_sum}
          q_{i,0}/z&=q_{i-1,1}+q_{i+1,1}+q_{i-1,0}+q_{i+1,0},&i&\geq 2,\\
          \label{eq_i=0_j>0_sum}
          q_{0,j}/z&=q_{1,j-1}+q_{1,j+1}+q_{0,j-1}+q_{0,j+1},&j&\geq 2,
     \end{align}
and
     \begin{align}
          \label{qq10}
          q_{1,0}/z&=q_{0,1}+q_{2,1}+q_{0,0}+q_{2,0},\\
          \label{qq01}
          q_{0,1}/z&=q_{1,0}+q_{1,2}+q_{0,0}+q_{0,2},\\
          \label{eq_i=0_j=0_sum}
          q_{0,0}/z&=1/z+q_{0,1}+q_{1,1}+q_{1,0}. & &
\end{align}
(For a reason that will be clear below, we want to discus separately the equations involving $q_{0,0}$, and we therefore do not merge \eqref{qq10} (resp.\ \eqref{qq01}) and \eqref{eq_i>0_j=0_sum} (resp.\ \eqref{eq_i=0_j>0_sum}).

\begin{lem}\label{lemm1}
Equations \eqref{eq_i>0_j>0_sum}--\eqref{eq_i=0_j=0_sum} have a unique solution in the form of formal power series.
\end{lem}
\begin{proof}
The generating function $q_{i,j}$ can be written as a Taylor series about $z=0$. Substituting these Taylor series into Equations~(\ref{eq_i>0_j>0_sum})--(\ref{eq_i=0_j=0_sum}), and equating coefficients of $z^k$, yields Equations~(\ref{eq_i>0_j>0})--(\ref{eq_i=0_j=0}). The latter obviously have a unique solution for the coefficients $q_{i,j,k}$, because these counting numbers can be determined recursively using $q_{0,0,0}=1$.
\end{proof}

In order to find the unique solution for $q_{i,j}$, we shall employ the compensation approach, which consists of three steps:
\begin{itemize}
\item Characterize all products $\alpha^{i}\beta^{j}$ for which the inner equations (\ref{eq_i>0_j>0_sum}) are satisfied, and construct linear combinations of these products, which in addition to being formal solutions to (\ref{eq_i>0_j>0_sum}),
also satisfy (\ref{eq_i>0_j=0_sum}) and (\ref{eq_i=0_j>0_sum}).
\item Prove that these solutions are formal power series.
\item Determine the complete unique solution $q_{i,j}$ by taking into account the boundary conditions \eqref{qq10}--\eqref{eq_i=0_j=0_sum}.
\end{itemize}

\subsection{Linear combinations of products}\label{product_forms}

\begin{figure}[!ht] \centering
 \includegraphics[width= .35\linewidth]{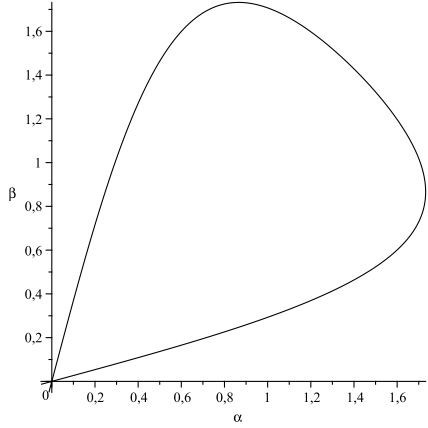}
 \caption{The curve (\ref{inner_condition}) in $\mathbb{R}_{+}^{2}$ for $z=1/4$. }
 \label{curve}
\end{figure}

Substituting the product $\alpha^i\beta^j$ into the inner equations (\ref{eq_i>0_j>0_sum}), and dividing by common powers, yields
     \begin{equation}
     \label{inner_condition}
          \alpha \beta/z=\beta^2 + \alpha^2+\alpha^2\beta^2.
     \end{equation}
Incidentally, note that \eqref{inner_condition} is nothing else but \eqref{dasas} (see \eqref{kern2}).  Hence, a product  $\alpha^i\beta^j$ is a solution of (\ref{eq_i>0_j>0_sum}) if and only if (\ref{inner_condition}) is satisfied, and any linear combination of such products will satisfy (\ref{eq_i>0_j>0_sum}).
Figure \ref{curve} depicts the curve (\ref{inner_condition}) in $\mathbb{R}_{+}^{2}$ for $z=1/4$.

Now we construct a linear combination of the products introduced
above, which will give a formal solution to the equations (\ref{eq_i>0_j>0_sum})--(\ref{eq_i=0_j>0_sum}).
The first term of this combination, say $\alpha_0^i\beta_0^j$, we require to satisfy both
(\ref{eq_i>0_j>0_sum}) and (\ref{eq_i>0_j=0_sum}).
In other words, the pair $(\alpha_0,\beta_0)$ has to satisfy (\ref{inner_condition})
as well as
     \begin{equation}
     \label{horizontal_condition}
          \alpha \beta/z=\beta^2+\alpha^2\beta^2+\beta+\alpha^2\beta.
     \end{equation}
     The motivation to start with a term satisfying both (\ref{inner_condition}) and \eqref{horizontal_condition} will be explained at the end of this subsection, see Remark \ref{rempair}.
\begin{lem}
\label{lem_a_b}
There exists a unique pair $(\alpha_0,\beta_0)$, with
     \begin{equation}
     \label{def_alpha_0}
          \alpha_0=\frac{1-\sqrt{1-8z^2}}{4z},
          \quad \beta_0=\frac{\alpha_0^2}{1+\alpha_0^2},
     \end{equation}
 of formal power series in $z$ such that the
product $\alpha_0^i\beta_0^j$ satisfies both \eqref{eq_i>0_j>0_sum} and \eqref{eq_i>0_j=0_sum}.
\end{lem}

\begin{proof}
Equation \eqref{horizontal_condition} yields
     \begin{equation}
     \label{intermediate_eq}
          \alpha/z=\beta+\alpha^2\beta+1+\alpha^2=(1+\beta)(1+\alpha^2).
     \end{equation}
Subtracting \eqref{horizontal_condition} from \eqref{inner_condition} gives
     $
          \beta(1+\alpha^{2})=\alpha^2,
     $
so that  $(1+\beta)(1+\alpha^2)=1+2\alpha^2$, and together
with (\ref{intermediate_eq}) this gives
$\alpha/z=1+2\alpha^2$. The last equation has the two solutions $(1\pm \sqrt{1-8z^2})/(4z)$, and  only the solution with the $-$ sign is a formal power series in $z$.
\end{proof}

The function $\alpha_0^i\beta_0^j$ thus satisfies (\ref{eq_i>0_j>0_sum})
and (\ref{eq_i>0_j=0_sum}), but it fails to satisfy
(\ref{eq_i=0_j>0_sum}). Indeed, if $\alpha^i\beta^j$
is a solution to (\ref{eq_i=0_j>0_sum}), then $(\alpha,\beta)$ should satisfy
     \begin{equation}
     \label{vertical_condition}
          \alpha \beta/z=\alpha^2+\alpha^2\beta^2+\alpha+\alpha\beta^2,
     \end{equation}
     and  $(\alpha_0,\beta_0)$ is certainly not a solution to (\ref{vertical_condition}).

\begin{figure}[!ht] \centering
 \includegraphics[width= .35\linewidth]{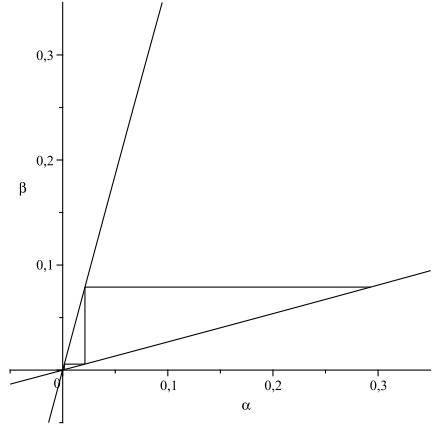}
 \caption{The curve (\ref{inner_condition}) near $0$ for $z=1/4$,  and the construction
 of the products.}
 \label{fig_construction}
\end{figure}

Now we start adding compensation terms. We consider
$c_0\alpha_0^i\beta_0^j+d_1\alpha^i\beta^j$, where $(\alpha,\beta)$
satisfies (\ref{inner_condition}) and is such that $c_0\alpha_0^i\beta_0^j+d_1\alpha^i\beta^j$
satisfies (\ref{eq_i=0_j>0_sum}).
The identity (\ref{eq_i=0_j>0_sum}), which has to be true for any $j\geq 2$, forces us to take
$\beta=\beta_0$. Then, thanks to (\ref{inner_condition}), we see
that $\alpha$ is the conjugate root of $\alpha_0$ in
(\ref{inner_condition}). Denote this new
root by $\alpha_1$.
Hence, $\alpha_0$ and $\alpha_1$ are the two roots of Equation (\ref{inner_condition}),
where $\beta$ is replaced by $\beta_0$. In other words, $\alpha_0$ and $\alpha_1$
satisfy
     \begin{equation}
          (1+\beta_0^2)\alpha^2-(\beta_0/z)\alpha+\beta_0^2=0.
     \end{equation}
Due to the root-coefficient relationships, we obtain
     \begin{equation}
     \label{before_conclusion}
          \frac{\beta_0/z}{1+\beta_0^2}=\alpha_0+\alpha_1.
     \end{equation}
We now determine $d_1$ in terms of $c_0$. Note that
$c_0\alpha_0^i\beta_0^j+d_1\alpha_1^i\beta_0^j$
is a solution to (\ref{eq_i=0_j>0_sum}) if and only if
     \begin{equation}
          (\beta_0/z)(c_0+d_1)=(\alpha_0 c_0+\alpha_1 d_1)(1+\beta_0^2)+
          (c_0+d_1)(1+\beta_0^2).
     \end{equation}
The latter identity can be rewritten as
     \begin{equation}
     \label{before_simplif}
          \frac{\beta_0/z}{1+\beta_0^2}(c_0+d_1)=(\alpha_0 c_0+\alpha_1 d_1)+
          (c_0+d_1),
     \end{equation}
and thus, using  \eqref{before_conclusion},
     \begin{equation}
     \label{simplif}
          d_1=-\frac{1-\alpha_1}{1-\alpha_0}c_0.
     \end{equation}
The series $c_0\alpha_0^i\beta_0^j+d_1\alpha_1^i\beta_0^j$ after one compensation step satisfies
 \eqref{eq_i>0_j>0_sum}
         and  \eqref{eq_i=0_j>0_sum} for the interior and the vertical boundary. However, the compensation term $d_1\alpha_1^i\beta_0^j$ has generated a new error at the horizontal boundary. To compensate for this, we must add another compensation term, and so on. In this way,
the compensation approach can be continued, which eventually leads to
\begin{equation}
\label{after_compensation}
x_{i,j}=\rlap{$\underbrace{\phantom{c_0\alpha_0^i\beta_0^j+
                                    d_1\alpha_1^i\beta_0^j}}_{\rm V}$}
\overbrace{c_0\alpha_0^i\beta_0^j}^{\rm H}+
        \rlap{$\overbrace{\phantom{d_1\alpha_1^i\beta_0^j+
                                   c_1\alpha_1^i\beta_1^j}}^{\rm H}$}
d_1\alpha_1^i\beta_0^j+
        \rlap{$\underbrace{\phantom{c_1\alpha_1^i\beta_1^j+
                                    d_2\alpha_2^1\beta_1^j}}_{\rm V}$}
c_1\alpha_1^i\beta_1^j+
\overbrace{d_2\alpha_2^i\beta_1^j+c_2\alpha_2^i\beta_2^j}^{\rm H}+\cdots
\end{equation}
where, by construction and (\ref{inner_condition}), for all $k\geq 0$ we have
     \begin{equation}
     \label{properties_alpha_beta}
          \beta_k=f(\alpha_k),
          \quad \alpha_{k+1}=f(\beta_k),
          \quad f(t)=\frac{1-\sqrt{1-4z^{2}(1+t^2)}}{2z(1+t^2)}t,
     \end{equation}
     and the function $f(t)$ follows from    \eqref{inner_condition}.
The construction is such that each term in \eqref{after_compensation} satisfies the equations in the
interior of the state space, the sum of two terms with the same $\alpha$ satisfies the horizontal
boundary conditions (H) and the sum of two terms with the same $\beta$ satisfies the vertical
boundary conditions (V). Hence $x_{i,j}$ satisfies all  equations  \eqref{eq_i>0_j>0_sum}--\eqref{eq_i=0_j>0_sum}.
See also Figure \ref{fig_construction}. Furthermore,
     \begin{align}
     \label{recurrence_coefficients}
          d_{k+1}&=-\frac{1-\alpha_{k+1}}{1-\alpha_k}c_k, \quad
          c_{k+1}=-\frac{1-\beta_{k+1}}{1-\beta_k}d_{k+1},
          \quad  k\geq 0.
     \end{align}
An easy calculation starting from (\ref{recurrence_coefficients}) yields
     \begin{equation}
          d_{k+1}=-\frac{(1-\alpha_{k+1})(1-\beta_k)}{(1-\alpha_0)(1-\beta_0)}c_0,
          \quad c_{k+1}=\frac{(1-\beta_{k+1})(1-\alpha_{k+1})}{(1-\alpha_0)(1-\beta_0)}c_0,
          \quad  k\geq 0,
     \end{equation}
so that choosing (arbitrarily) $c_0=(1-\alpha_0)(1-\beta_0)$ finally gives
     \begin{equation}
     \label{formal_solution_1}
          x_{i,j}=\sum_{k=0}^{\infty}(1-\beta_k)\beta_k^j
          [(1-\alpha_k)\alpha_k^i-(1-\alpha_{k+1})\alpha_{k+1}^i].
     \end{equation}

By symmetry, we also obtain that $x_{j,i}$ is a formal solution to \eqref{eq_i>0_j>0_sum}--\eqref{eq_i=0_j>0_sum}. This leaves to consider the solution for the states $(0,1)$, $(0,0)$ and $(1,0)$, i.e., the boundary conditions \eqref{qq10}--\eqref{eq_i=0_j=0_sum}.

\subsection{Formal power series}

\begin{prop}
\label{prop_control}
The sequences $\{\alpha_k\}_{k\geq 0}$ and $\{\beta_k\}_{k\geq 0}$ appearing in \eqref{after_compensation}--\eqref{properties_alpha_beta}
satisfy the following property: for all $k\geq 0$, $\alpha_k=z^{2k+1}\widehat{\alpha}_k(z)$  and
$\beta_k=z^{2k+2}\widehat{\beta}_k(z)$, where $\widehat{\alpha}_k$ and $\widehat{\beta}_k$
are formal power series such that $\widehat{\alpha}_k(0)=\widehat{\beta}_k(0)=1$.
\end{prop}

\begin{proof}
For $f$ defined in (\ref{properties_alpha_beta}), all $p\geq 0$, and all real numbers $s_1,s_2,\ldots$,
     \begin{equation}
          f(z^p[1+s_{1}z+s_{2}z^{2}+\cdots ])=z^{p+1}[1+s_1z+(1+s_2)z^2+\cdots ].
     \end{equation}
Because (\ref{def_alpha_0}) yields
$\alpha_0=z+2z^3+\cdots$, the proof is completed via (\ref{properties_alpha_beta}).
\end{proof}

The following result is an immediate consequence of Proposition \ref{prop_control}.
\begin{cor}
\label{convergence_everywhere}
For  $i,j \geq 0$, $x_{i,j}$ and $x_{j,i}$ defined in \eqref{formal_solution_1} are formal power series in $z$.
\end{cor}

\begin{rem}\label{rempair}
The approach outlined in Subsection \ref{product_forms} is initialized with a term satisfying both
\eqref{inner_condition} and \eqref{horizontal_condition}. Alternatively, we could start with an arbitrarily chosen term with power series $\alpha_0$ and $\beta_0$ satisfying \eqref{inner_condition} only. This term would violate  \eqref{horizontal_condition} as well as \eqref{vertical_condition}, and therefore generate two sequences of terms, one starting with compensation of \eqref{horizontal_condition} and the other with \eqref{vertical_condition}. It is readily seen that in one of the two sequences, some of the parameters $\alpha$ or $\beta$ will not be power series. Indeed, the two sequences of terms would be generated as in \eqref{properties_alpha_beta}, with instead of $f$ the function 
\begin{equation}
     f_\pm(t)=\frac{1\pm\sqrt{1-4z^{2}(1+t^2)}}{2z(1+t^2)}t.
\end{equation}
In Subsection \ref{product_forms}, we already described the construction of the sequence with the choice $f_-=f$. For $f=f_+$, notice that for all $p\geq 0$, and all real numbers $s_1,\ldots$, we have
     \begin{equation}
          f_+(z^p[1+s_{1}z+\cdots ])=z^{p-1}[1+s_1z+\cdots ].
     \end{equation}
Hence, the resulting solution constructed as in \eqref{properties_alpha_beta} would fail to be a formal power series, cf.~Equation \eqref{formal_solution_1} and Subsection \ref{subsec:determining}.
\end{rem}

\subsection{Determining the unique solution}
\label{subsec:determining}

We now determine the generating functions $q_{i,j}$. Define
     \begin{equation}
     \label{defhat}
          \widehat{x}_{i,j}=x_{i,j}+x_{j,i},
     \end{equation}
so that in particular
\begin{equation}
    \label{ex_x00}
         \widehat{x}_{0,0}=2\sum_{k=0}^{\infty}(1-\beta_k)(\alpha_{k+1}-\alpha_k).
\end{equation}

\begin{prop}\label{sds}
The expressions
\begin{equation}
     \label{value_q00}
          q_{0,0}=\frac{1+\widehat{x}_{0,0}}{1-2z+z\widehat{x}_{0,0}}
     \end{equation}
and
     \begin{equation}
     \label{value_qij}
          q_{i,j}=c \widehat{x}_{i,j},\quad i+j>0,
     \end{equation}
     with
        \begin{equation}
    \label{value_c}
          c=\frac{1}{1-2z+z\widehat{x}_{0,0}},
     \end{equation}
are the unique solutions to Equations {\rm(\ref{eq_i>0_j>0_sum})--(\ref{eq_i=0_j=0_sum})} in the form of formal power series.
\end{prop}

\begin{proof}
We shall look for a solution satisfying, for $i+j>0$,
\begin{equation}
q_{i,j} = c x_{i,j}+\widetilde{c} x_{j,i}.
\end{equation}
By symmetry we must have $c=\widetilde{c}$.
Such combinations have been shown to satisfy
 \eqref{eq_i>0_j>0_sum}--\eqref{eq_i=0_j>0_sum}, since they do not involve $q_{0,0}$. We should still determine the solution for which also the boundary conditions \eqref{qq10} and \eqref{eq_i=0_j=0_sum} (or by symmetry \eqref{qq01}) are satisfied. This can be achieved by choosing $c$ and $q_{0,0}$ appropriately.

By using (\ref{eq_i=0_j=0_sum})
and (\ref{qq10}) (or  (\ref{eq_i=0_j=0_sum}) and \eqref{qq01}), we obtain that
     \begin{equation}
     \label{value_cc}
          c=\frac{z}{\widehat{x}_{1,0}-z[\widehat{x}_{0,1}+\widehat{x}_{2,1}+\widehat{x}_{2,0}]-
          z^2[\widehat{x}_{1,0}+\widehat{x}_{0,1}+\widehat{x}_{1,1}]}
     \end{equation}
and
     \begin{equation}
     \label{value_q000}
          q_{0,0}=1+zc[\widehat{x}_{0,1}+\widehat{x}_{1,1}+\widehat{x}_{1,0}].
     \end{equation}
In addition, tedious calculations starting from \eqref{formal_solution_1} give
     \begin{equation}
     \label{tr}
          \widehat{x}_{1,0}/z-[\widehat{x}_{0,1}
          +\widehat{x}_{2,1}+\widehat{x}_{2,0}+\widehat{x}_{0,0}]=1,
          \quad (1/z-1)\widehat{x}_{0,0}+2-[\widehat{x}_{1,0}+
          \widehat{x}_{0,1}+\widehat{x}_{1,1}]=0.
     \end{equation}
Thanks to the equations in \eqref{tr}, the denominator of  \eqref{value_cc} can be written as $z[1-2z+z\widehat{x}_{0,0}]$.
The functions defined in \eqref{value_q00} and \eqref{value_qij} are the solutions to Equations {\rm(\ref{eq_i>0_j>0_sum})--(\ref{eq_i=0_j=0_sum})} by construction. Finally, the uniqueness follows from Lemma \ref{lemm1}.
\end{proof}


\begin{prop}
\label{prop_explicit_expression}
     \begin{align}
          Q(x,y;z)&=\sum_{i,j=0}^{\infty} q_{i,j}x^i y^j\label{fieq}\\
                  &=c+ c\sum_{k=0}^{\infty} \left(\frac{1-\beta_k}{1-\beta_k y}
          \left[\frac{1-\alpha_k}{1-\alpha_k x}     -
          \frac{1-\alpha_{k+1}}{1-\alpha_{k+1} x}\right]
          +
          \frac{1-\beta_k}{1-\beta_k x}
          \left[\frac{1-\alpha_k}{1-\alpha_k y}     -
          \frac{1-\alpha_{k+1}}{1-\alpha_{k+1} y}\right]\right),\label{seeq}
     \end{align}
with $q_{i,j}$ as in Proposition {\rm \ref{sds}}, $c$ as in \eqref{value_c}, and $\alpha_k$ and $\beta_k$ as in \eqref{def_alpha_0} and \eqref{properties_alpha_beta}, respectively.
Furthermore,
\begin{equation}
\label{eq:simplifed-functions}
Q(0,0;z) = c(1+\widehat x_{0,0}),
\quad
Q(1,0;z)=Q(0,1;z)= c(1-\alpha_0(z)),
\quad
Q(1,1;z)=c.
\end{equation}
\end{prop}


\begin{proof}
This is immediate from (\ref{formal_solution_1}), \eqref{ex_x00}--\eqref{value_c}.
\end{proof}

\begin{prop}
\label{infinite_poles}
For all $0<|z|<1/3$, the functions $\sum_{i,j=0}^{\infty} q_{i,j}x^i y^j$, $\sum_{i=0}^{\infty} q_{i,0}x^i$ and $\sum_{j=0}^{\infty} q_{0,j}y^j$ have infinitely many poles.
\end{prop}

\begin{proof}
This is a direct consequence of Propositions \ref{prop_control} and \ref{prop_explicit_expression}.
\end{proof}

It follows immediately from Proposition \ref{infinite_poles} that both $\sum_{i=0}^{\infty} q_{i,0}x^i$ and $\sum_{j=0}^{\infty} q_{0,j}y^j$ are nonholonomic. The holonomy being maintained after specialization to a variable (see \cite{FLAJ}), we reach the conclusion that the trivariate generating function $Q(x,y;z)$ is nonholonomic as well.

%

\subsection{Retrieving coefficients}
\label{Retrieving_coefficients}
Now that we have an explicit expression for $q_{i,j}$, we briefly present an efficient procedure for calculating its coefficients $q_{i,j,k}$.

To compute $q_{i,j}$ we need to calculate in principle an infinite series, or to solve the recurrence relations. However, we shall prove that if we are only interested in a finite number of coefficients $q_{i,j,k}$, then it is enough to take into account a finite number of $\alpha_k$ and $\beta_k$. For $k\geq 0$, denote by $k\vee 1$ the maximum of $k$ and $1$. Define
\begin{equation}
\label{tresh}
N_p^{i,j}=1+\left\lfloor\tfrac{1}{4}\max\{p-(i\vee1 +2(j\vee1)),p-(2(i\vee1)+j\vee1)\}\right\rfloor.
\end{equation}
\begin{prop}
\label{prop_enough}
For any $i,j\geq 0$, the first $p$ coefficients of $q_{i,j}$ only require the series expansions of order $p$ of
     $\alpha_0,\beta_0,\ldots \alpha_{N_p^{i,j}},\beta_{N_p^{i,j}}$.
\end{prop}

\begin{proof}
In order to obtain the series expansion of $q_{0,0}$ of order $p$,
it is enough to know the series expansion of $\widehat{x}_{0,0}$
of order $p$, see \eqref{value_q00}.
From  \eqref{formal_solution_1} we obtain
     \begin{equation}
     \label{nfx}
          \widehat{x}_{0,0}=2\sum_{k=0}^{\infty}(1-\beta_k)(\alpha_{k+1}-\alpha_k)
          =-2\alpha_0+2\sum_{k=0}^{\infty}\beta_k(\alpha_{k}-\alpha_{k+1}).
     \end{equation}
With Proposition \ref{prop_control}, $\beta_k\alpha_k=O(z^{4k+3})$ and
$\beta_k\alpha_{k+1}=O(z^{4k+5})$, so that in order to obtain the
series expansion of $\widehat{x}_{0,0}$ of order $p$, it is enough
to consider in \eqref{nfx} the values of $k$ such that $4k+3\leq p$.
In other words, we have to deal with $\alpha_0,\beta_0,\ldots ,\alpha_k,\beta_k,\alpha_{k+1}$
for $4k+3\leq p$.
Since $N_{p}^{0,0}=1+\lfloor (p-3)/4\rfloor$, see \eqref{tresh},
Proposition \ref{prop_enough} is shown for $i=j=0$.
The proof for other values of $i$ and $j$ is similar and therefore omitted.
\end{proof}

As an application of Proposition \ref{prop_enough}, let us find
the numbers  $q_{0,0,k}$ for $k\in\{0,\ldots ,10\}$.
Taking $i=j=0$, we have $N_{10}^{0,0}=2$; we then calculate
the series expansions of order $10$ of $\alpha_0,\beta_0,
\alpha_1,\beta_1,\alpha_2,\beta_2$.
After using (\ref{value_q00})
and \eqref{nfx}, we obtain the numbers:
     \begin{equation}
          1,0,2,2,10,16,64,126,454,1004,3404.
     \end{equation}

\section{Asymptotic analysis}
\label{Asymptotic_analysis}

In this section we derive an asymptotic expression for the coefficients $q_{i,j,k}$ for large values of $k$, using the technique of \textit{singularity analysis} (see Flajolet and Sedgewick \cite{FLAJ} for an elaborate exposition). This requires the investigation of the function $q_{i,j}=q_{i,j}(z)$ near its dominant singularity (closest to the origin) in the $z$-plane. First, notice that
     \begin{equation}
     \label{finitesum}
          \sum_{i,j=0}^{\infty}q_{i,j,k}\leq 3^k,
     \end{equation}
because there are at most $3^k$ paths of length $k$. This clearly implies that $q_{i,j,k}\leq 3^k$, so that $q_{i,j}$ is analytic at least for $|z|<1/3$, and the singularities $z$ of $q_{i,j}$ must satisfy $|z|\geq1/3$.

In fact, the singularities of $q_{i,j}$ are given by the singularities of $\widehat{x}_{i,j}$ and the zeros of the denominator of $c$ in \eqref{value_c}, see \eqref{value_q00} and \eqref{value_qij} . In what follows, we first show in Subsection \ref{sub:conv} that $\widehat{x}_{i,j}$ converges in the domain $\vert z\vert \leq 1/\sqrt8$. We then find in Subsection \ref{sub:fs} the first singularity of $q_{i,j}$. Finally, we state and prove our asymptotic results in Subsection \ref{sub:asre}.

\subsection{Convergence of the solutions}
\label{sub:conv}
Since $\widehat{x}_{i,j}$ is constructed from the functions $\alpha_0,\beta_0,\alpha_1,\beta_1,\ldots$, and since all these functions follow from the iterative scheme \eqref{properties_alpha_beta}, i.e., from $\alpha_0$, it is clear that the singularities of $\widehat{x}_{i,j}$ are in fact the singularities of $\alpha_0$, namely (see \eqref{def_alpha_0})
    \begin{equation}
     z=\pm \frac{1}{\sqrt{8}}.
     \end{equation}
Further, we have that:

\begin{prop}
\label{prop:analytic-properties}
The sequences $\{\alpha_k\}_{k\geq 0}$ and $\{\beta_k\}_{k\geq 0}$ appearing in \eqref{after_compensation}--\eqref{properties_alpha_beta} satisfy the following properties:
\begin{enumerate}[label={\rm (\roman{*})},ref={\rm (\roman{*})}]
     \item \label{entrelace}
           For all $|z|\leq 1/\sqrt8$, we have $1/\sqrt 2\geq |\alpha_0|>|\beta_0|>|\alpha_1|>|\beta_1|>\cdots $.
     \item \label{exponential_control}
           For all $|z|\leq 1/\sqrt8$ and all $k\geq 0$, we have $|\alpha_k|\leq 1/\sqrt 2^{2k+1}$ and $|\beta_k|\leq 1/\sqrt 2^{2k+2}$.
\end{enumerate}
\end{prop}
Before starting  with the proof of Proposition \ref{prop:analytic-properties}, let us recall Rouch\'e's theorem (see e.g.~\cite[Chapter 4]{FLAJ}).
\begin{thm}[Rouch\'e's theorem]
\label{RTT}
Let the functions $F$ and $G$ be analytic in a simply connected domain of $\mathbb C$ containing in its interior the closed simple curve $\gamma$ . Assume that $F$ and $G$ satisfy $\vert F(r )-G(r )\vert < \vert G(r )\vert$ for $r$ on the curve $\gamma$ . Then $F$ and $G$ have the same number of zeros inside the interior domain delimited by $\gamma$ .
\end{thm}

\begin{proof}[Proof of Proposition \ref{prop:analytic-properties}]
To show \ref{entrelace} and \ref{exponential_control}, it is enough, thanks to \eqref{properties_alpha_beta} and since $\alpha_0 = f(1)$, to prove that for all $|t|\leq 1$ and $|z|\leq 1/\sqrt8$,
     \begin{equation}
     \label{properties_f}
          |f(t)|\leq \frac{|t|}{\sqrt2}.
     \end{equation}
In order to show (\ref{properties_f}), we
first note that $f$ satisfies the algebraic equation
$(1+t^2)f(t)^2-tf(t)/z+t^2=0$, see (\ref{inner_condition}) and \eqref{properties_alpha_beta}.
As a consequence, the function
$f(t)/t$ is such that (with $r=f(t)/t$)
     \begin{equation}
     \label{algebraic_g}
         (1+t^2)r^2-r/z+1=0.
     \end{equation}
Let us now fix $\vert t\vert \leq 1$, $\vert z\vert<1/\sqrt 8$, and introduce the functions
     \begin{equation}
          F(r ) = (1+t^2)r^2-r/z+1,
          \quad
          G(r ) = -r/z+1.
     \end{equation}
We have
     \begin{equation}
          \vert F(r ) -G(r ) \vert = \vert (1+t^2) r^2\vert\leq 2\vert r\vert^2,
          \quad
          \vert G(r )\vert = \vert -r/z+1\vert \geq \vert r\vert/\vert z\vert -1.
     \end{equation}
Further, a straightforward computation gives that for any
     \begin{equation}
          \frac{1-\sqrt{1-8\vert z\vert^2}}{4\vert z\vert} < r_z < \frac{1+\sqrt{1-8\vert z\vert^2}}{4\vert z\vert}
     \end{equation}
we have, for all $\vert r\vert =r_z$,
     \begin{equation}
          2\vert r\vert^2 < \vert r\vert/\vert z\vert -1.
     \end{equation}
It then follows from Rouch\'{e}'s theorem, applied to the circle with radius $r_z$, that
\eqref{algebraic_g} has one unique solution $r=f(t)/t$ with $|r|< r_z$. By a continuity argument, we deduce that for all $\vert t\vert\leq 1$ and all $\vert z\vert <1/\sqrt8$,
     \begin{equation}
          \frac{\vert f(t)\vert}{\vert t\vert}\leq \frac{1-\sqrt{1-8\vert z\vert^2}}{4\vert z\vert}.
     \end{equation}
Equation \eqref{properties_f} and Proposition \ref{prop:analytic-properties} follow immediately.
\end{proof}

\begin{cor}
\label{cor:ce}
For  $i,j \geq 0$, $x_{i,j}$ and $x_{j,i}$ defined in \eqref{formal_solution_1} are convergent for $|z|\leq 1/\sqrt8$. Further, for each $|z|\leq1/\sqrt8$,
\begin{equation}
         \sum_{i,j= 0}^{\infty}|x_{i,j}(z)|< \infty.
     \end{equation}
\end{cor}

\begin{proof}
The first statement in Corollary \ref{cor:ce} immediately follows from Proposition \ref{prop:analytic-properties} and from the expression for $x_{i,j}$. Next, use \eqref{formal_solution_1} to get that for $i+j>0$,
     \begin{equation}
          |x_{i,j}|\leq \sum_{k=0}^{\infty}|1-\beta_k||\beta_k|^j
          \left[|1-\alpha_k||\alpha_k|^i+|1-\alpha_{k+1}||\alpha_{k+1}|^i\right].
     \end{equation}
Using then that $|1-\beta_k|$, $|1-\alpha_k|$ and $|1-\alpha_{k+1}|$ are bounded, see Proposition \ref{prop:analytic-properties}\ref{entrelace}, and applying Proposition \ref{prop:analytic-properties}\ref{exponential_control} gives that there is some constant constant $C>0$ such that for all $i+j >0$,
     \begin{equation}
     \label{upper_bound_x}
          |x_{i,j}|\leq C\sum_{k=0}^{\infty}
          1/\sqrt2^{(2k+1)i+(2k+2)j}.
     \end{equation}
Likewise, we can prove that for some constant $C>0$,
     \begin{equation}
          |x_{0,0}|\leq C\sum_{k=0}^{\infty}1/\sqrt2^{(2k+1)}.
     \end{equation}
The second claim in Corollary \ref{cor:ce} immediately follows.
\end{proof}

\subsection{Dominant singularity of $q_{i,j}$}
\label{sub:fs}

A consequence of the results in Subsection \ref{sub:conv} is that $q_{i,j}$ is analytic for $\vert z\vert <1/\sqrt8$, except possibly at points $z$ such that the denominator of $c$ in \eqref{value_c} equals $0$. Denote this denominator by
      \begin{equation}
      \label{def_h}
           h(z)=1-2z+z\widehat{x}_{0,0},
      \end{equation}
and let $\rho$ denote the dominant singularity of $c$.
\begin{lem}
\label{def_rho}
The series $h(z)$ has a unique root in the interval $(1/3,1/\sqrt{8})$, which is the dominant singularity $\rho$ of $c$.
\end{lem}
The proof of Lemma \ref{def_rho} is presented in Appendix \ref{fa}. Furthermore, in Corollary \ref{approx_rho} we shall prove that $\rho\in[0.34499975,0.34499976]$

\begin{rem}
It seems hard to find a closed-form expression for $\rho$ or to decide whether $\rho$ is algebraic or not. From this point of view, the situation is quite different from that for the $74$ (nonsingular) models of walks discussed in Section \ref{sec:intro}. Indeed, for all these $74$ models, the analogue of $\rho$ (which eventually will be the first positive singularity of $Q(0,0;z)$, $Q(1,0;z)$, $Q(0,1;z)$ and $Q(1,1;z)$, see below) is always algebraic (of some degree between $1$ and $7$), see \cite{FRAofA}.
\end{rem}

\subsection{Asymptotic results}
\label{sub:asre}
Here are our main results on the asymptotic behavior of large counting numbers:
\begin{prop}
\label{prop_asymp_ij}
The asymptotics of the numbers of walks $q_{i,j,k}$ as $k\rightarrow \infty$ is
      \begin{equation}
      \label{main_result_asymp}
           q_{i,j,k}\sim C_{i,j}\rho^{-k},
      \end{equation}
with $\rho$ as defined in Lemma \ref{def_rho},
      \begin{equation}
C_{0,0}=\frac{3\rho-1}
                  {-\rho^2 h'(\rho)},
      \end{equation}
and
      \begin{equation}
           C_{i,j}=\frac{\widehat{x}_{i,j}(\rho)}{-\rho h'(\rho)},\quad i+j\geq 1.
      \end{equation}
\end{prop}

\begin{proof}
Consider first the case $i=j=0$.
Thanks to (\ref{value_q00}), \eqref{def_h} and Lemma \ref{def_rho},
we obtain that $q_{0,0}$ has a pole at $\rho$, and is holomorphic within the domain
     \begin{equation}
     \label{def_domain}
          \{z\in\mathbb{C}: |z|<(1+\epsilon)\rho\}\setminus [\rho,(1+\epsilon)\rho)
     \end{equation}
for any $\epsilon>0$ small enough. Moreover, the pole of $q_{0,0}$ at $\rho$ is of order one,
as we show now. For this it is sufficient to prove
that $1+\widehat{x}_{0,0}(\rho)\neq 0$ and that $h'(\rho)\neq 0$,
see again (\ref{value_q00}) and \eqref{def_h}.
In this respect, note that with \eqref{def_h} and Lemma \ref{def_rho} we have
$h(\rho)=1-2\rho+\rho\widehat{x}_{0,0}(\rho)=0$, so that
$1+\widehat{x}_{0,0}(\rho)=(3\rho-1)/\rho$,
which is positive by Lemma \ref{def_rho}.
On the other hand, it is a consequence of Lemma \ref{approx_mul}
that $h'(\rho)\neq  0$.
In particular, the behavior of $q_{0,0}$ near $\rho$ is given by
     \begin{equation}
     \label{local_behavior}
          q_{0,0}=\frac{1+\widehat{x}_{0,0}(\rho)}
                  {h'(\rho)(z-\rho)[1+O(z-\rho)]}
                 =\frac{3\rho-1}
                  {-\rho^2 h'(\rho)(1-z/\rho)}[1+O(z-\rho)].
     \end{equation}
The holomorphy of $q_{0,0}$ within the domain \eqref{def_domain} and the behavior
\eqref{local_behavior} of $q_{0,0}$ near $\rho$ immediately give
the asymptotics \eqref{main_result_asymp} for $i=j=0$ (see e.g.\ \cite{FLAJ}). The proof for other values of
$i$ and $j$ is similar and therefore omitted.
\end{proof}

A result similar to Proposition \ref{prop_asymp_ij} also holds for the total number of walks of length $k$ (not necessarily ending in state $(0,0)$). Indeed, from formula \eqref{eq:simplifed-functions}, we obtain that
\begin{equation}
Q(1,0;z)=\frac{1-\alpha_0(z)}{h(z)},
\quad
Q(1,1;z)=\frac{1}{h(z)},
\end{equation}
with $h(z)$ as in \eqref{def_h}. Accordingly, the exponential growth rate of these generating functions is the same as that of $Q(0,0;z)$, namely $\rho$ (defined in Lemma \ref{def_rho}). An identical reasoning as in the proof of Proposition \ref{prop_asymp_ij} then gives the following result.

\begin{prop}
\label{prop_asymp_11}
The asymptotics of the total number of walks $\sum_{i,j=0}^\infty q_{i,j,k}$ as $k\rightarrow \infty$ is
      \begin{equation}
      \label{main_result_asymp_11}
           \sum_{i,j=0}^\infty q_{i,j,k}\sim \frac{1}
                  {-\rho h'(\rho)}\rho^{-k}.
      \end{equation}
The asymptotics of the total number of walks ending on the horizontal axis (or, equivalently, on the vertical axis) $\sum_{i=0}^\infty q_{i,0,k}$ as $k\rightarrow \infty$ is
      \begin{equation}
      \label{main_result_asymp_10}
           \sum_{i=0}^\infty q_{i,0,k}\sim \frac{1-\alpha_0(\rho)}{-\rho h'(\rho)}\rho^{-k}.
      \end{equation}
\end{prop}

\begin{rem}
A priori, the fact that
\begin{equation}
\label{eq:fourfonctions}
     Q(0,0;z),\quad Q(1,0;z),\quad Q(0,1;z),\quad Q(1,1;z)
\end{equation}
have the same first positive singularity (namely, $\rho$) is far from obvious. This appears to be related to the fact that the two coordinates of the drift vector of our $\mathcal S$ are negative (they both equal $-2$). Here, for any step set $\mathcal S$, the drift vector is given by
$
     (\sum_{-1\leq i\leq 1} \delta_{i,j}^\mathcal S,\sum_{-1\leq j\leq 1} \delta_{i,j}^\mathcal S),
$
where $\delta_{i,j}^\mathcal S=1$ if $(i,j)\in\mathcal S$, and $0$ otherwise. Indeed, the following phenomenon has been observed in \cite{FRAofA} for the $79$ models of walks mentioned in Section \ref{sec:intro}: if the two coordinates of the drift vector are negative, then the four functions \eqref{eq:fourfonctions} do have the same first positive singularity.
\end{rem}

Let us illustrate Proposition \ref{prop_asymp_ij} with the following table, obtained by approximating $\rho$ by $0.34499975$  and $C_{0,0}$ by $0.0531$, see Corollary \ref{approx_rho} and Lemma \ref{approx_mul} in Appendix~\ref{fa}.

\medskip

\begin{center}
\begin{tabular}{|c|c|c|c|}
\hline
   Value of {$k$} & Exact value of $q_{0,0,k}$ & Approximation of $q_{0,0,k}$ & Ratio\\
    & (see Section \ref{Retrieving_coefficients}) & (see Proposition
\ref{prop_asymp_ij}) & \\
   \hline
   $10$ &  $3.404\cdot 10^{3\phantom{1}}   $ & $2.222\cdot
10^{3\phantom{1}} $ & $0.653$  \\
   $20$ &  $1.106\cdot 10^{8\phantom{1}}   $ & $9.305\cdot
10^{7\phantom{1}} $ & $0.840$  \\
   $30$ &  $4.254\cdot 10^{12}$ & $3.895\cdot 10^{12}$ & $0.915$  \\
   $40$ &  $1.714\cdot 10^{17}$ & $1.630\cdot 10^{17}$ & $0.951$  \\
   $50$ &  $7.037\cdot 10^{21}$ & $6.825\cdot 10^{21}$ & $0.969$  \\
   $60$ &  $2.913\cdot 10^{26}$ & $2.857\cdot 10^{26}$ & $0.980$  \\
   $70$ &  $1.211\cdot 10^{31}$ & $1.196\cdot 10^{31}$ & $0.987$  \\
   $80$ &  $5.051\cdot 10^{35}$ & $5.007\cdot 10^{35}$ & $0.991$  \\
   $90$ &  $2.109\cdot 10^{40}$ & $2.096\cdot 10^{40}$ & $0.993$  \\
   $100$ & $8.814\cdot 10^{44}$ & $8.774\cdot 10^{44}$ & $0.995$  \\
   \hline
\end{tabular}

\end{center}

\section{Discussion}
\label{Discussion}

\subsection{A wider range of applicability}
The compensation approach for our counting problem leads to an exact expression for the generating function $Q(x,y;z)$.
A detailed exposition of the compensation  approach can be found in \cite{MR1138205,MR1080417,MR1241929}, in which it has been shown to work for two-dimensional random walks on the lattice of the first quadrant that obey the following conditions:
\begin{itemize}
\item Step size: Only transitions to neighboring states.
\item Forbidden steps: No transitions from interior states to the North, North-East, and East.
\item Homogeneity: The same transitions occur for all interior points, and similarly for all points on the horizontal boundary, and for all points on the vertical boundary.
\end{itemize}
Although the theory has been developed for the {\it bivariate} transform of stationary distributions of recurrent two-dimensional random walks, this paper shows that is can also be applied to obtain a series expression for the {\it trivariate} function $Q(x,y;z)$. The fact that $Q(x,y;z)$ has one additional variable does not seem to matter much. We therefore expect that the compensation approach will work for walks in the quarter plane that obey the three conditions mentioned above. This will be a topic for future research. Another topic is to see whether the above conditions can be relaxed, and in fact, in the next section we present an example of a walk, not satisfying these properties, yet amenable to the compensation approach.

\subsection{Another example}
\begin{figure}[!ht]
 \begin{picture}(10.00,75.00)
 \hspace{37mm}
 \includegraphics{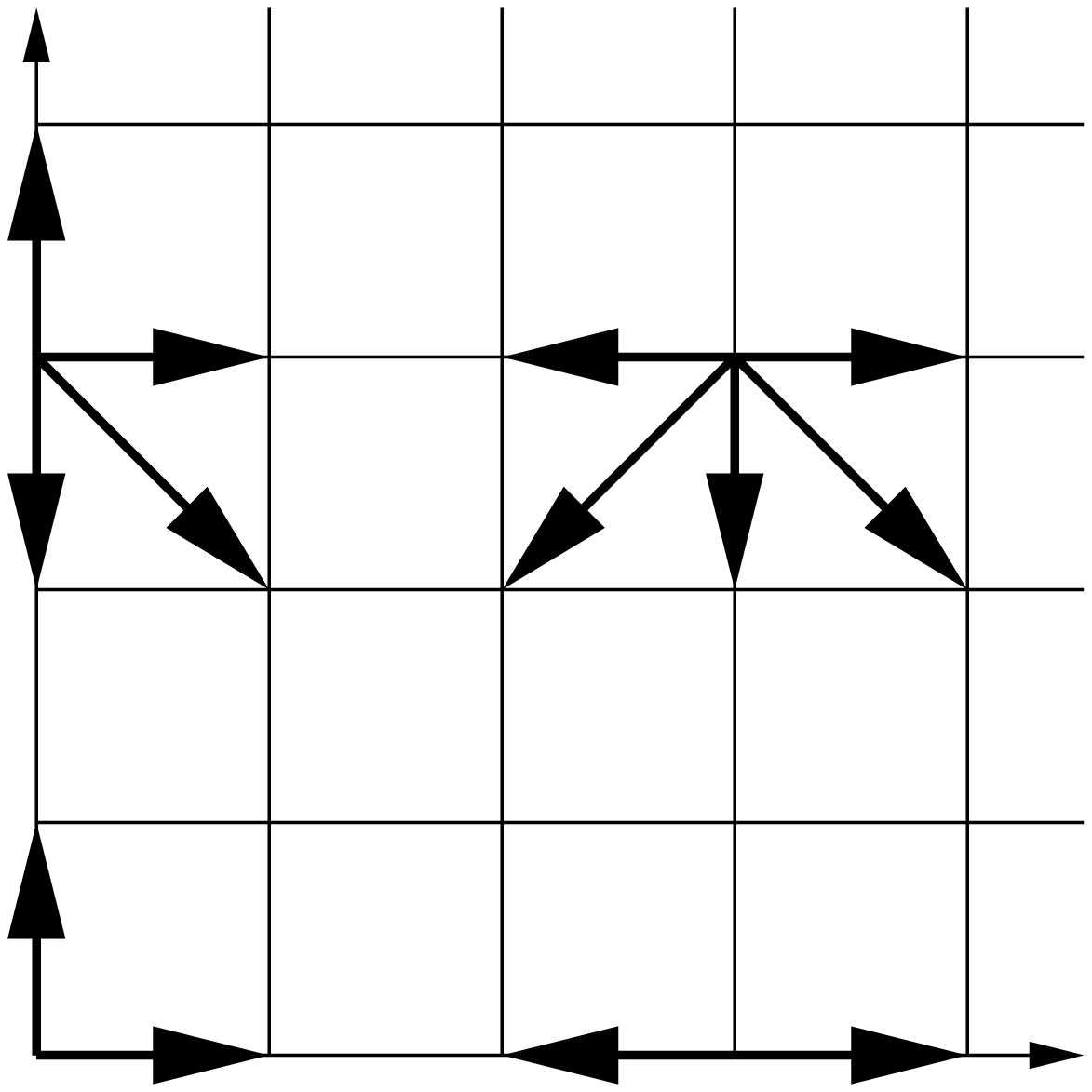}
 \hspace{60mm}
 \includegraphics{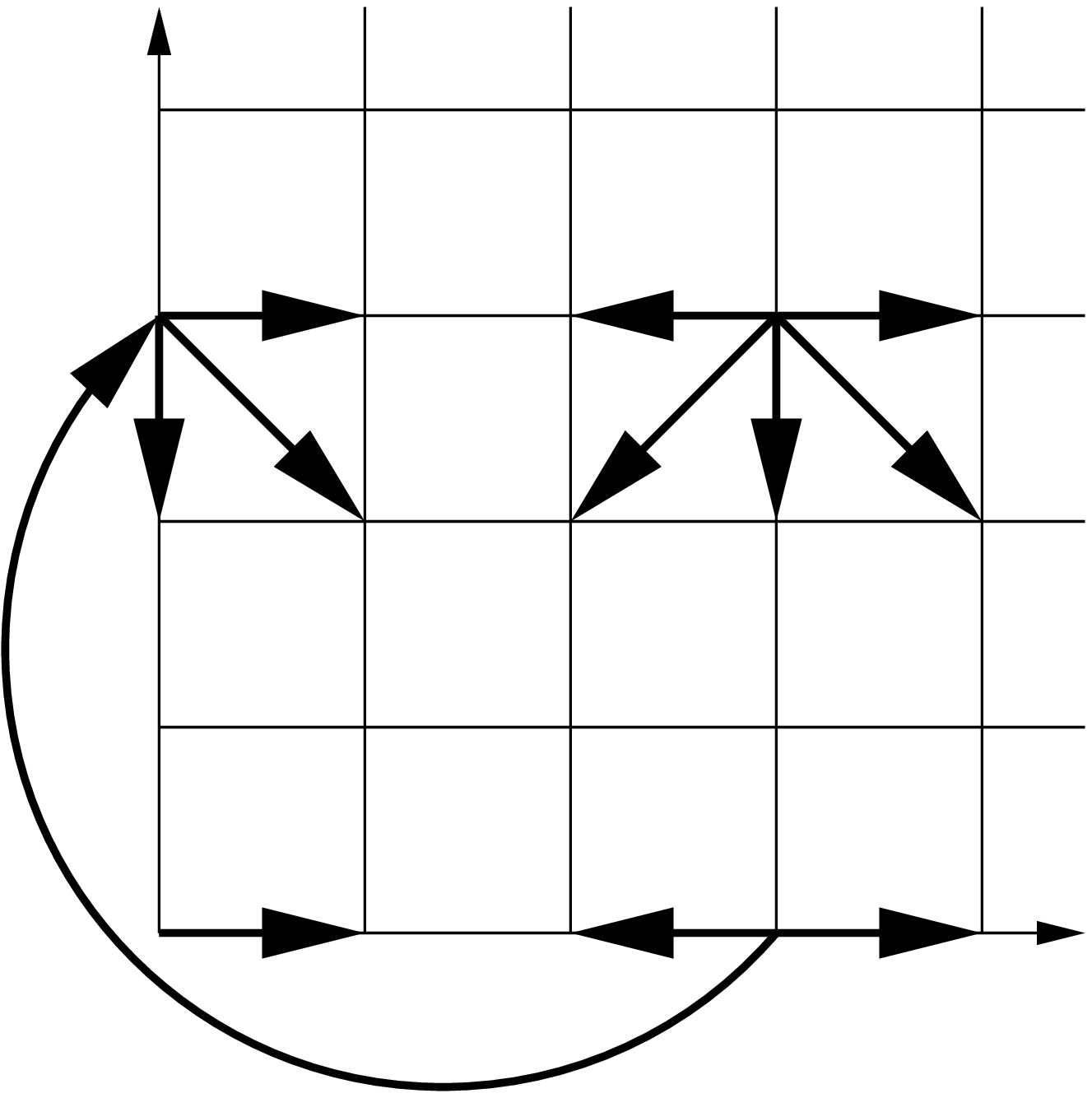}
 \end{picture}
 \caption{Two walks with different boundary behavior}
 \label{other_examples}
 \end{figure}

We consider the walk with $\mathcal{S}=\{(-1,0),(-1,-1),(0,-1),(1,-1),(1,0)\}$ in the
interior, $\mathcal{S}_H=\{(-1,0),(1,0)\}$ on the horizontal boundary, $\mathcal{S}_V=\{(0,1),(0,-1),(1,-1),(1,0)\}$ on the vertical boundary, and $\mathcal{S}_0=\{(0,1),(1,0)\}$; see the left display in Figure~\ref{other_examples}. This walk has a rather special behavior, as it can only move upwards on the vertical boundary. By simple enumeration, we get the following recursion relations:
\begin{align}
          \label{eq_i>0_j>0_sumS}
          q_{i,j}/z&=q_{i-1,j}+q_{i-1,j+1}+q_{i,j+1}+q_{i+1,j}+q_{i+1,j+1}, &i,j&>0,\\
          \label{eq_i>0_j=0_sumS}
          q_{i,0}/z&=q_{i-1,0}+q_{i-1,1}+q_{i,1}+q_{i+1,0}+q_{i+1,1}, &i&>0,\\
          \label{eq_i=0_j>0_sumS}
          q_{0,j}/z&=q_{0,j-1}+q_{0,j+1}+q_{1,j}+q_{1,j+1}, &j&>0,\\
          \label{eq_i=0_j=0_sumS}
          q_{0,0}/z&=1/z+q_{0,1}+q_{1,1}+q_{1,0}. & &
\end{align}

\begin{prop}
The unique solution to Equations \eqref{eq_i>0_j>0_sumS}--\eqref{eq_i=0_j=0_sumS} is given by $q_{i,j}=c_0\alpha_0^i\beta_0^j$, with
$\beta_0$ defined as the unique formal power series in $z$ and solving
\begin{equation}
\label{rel:bet}
\beta/z=1+\beta^2+\beta^2(1+\beta)^2,
\end{equation}
$\alpha_0=\beta_0(1+\beta_0)$ and
$
c_0=\beta_0/z
$.
Therefore, the GF is algebraic {\rm(}and even rational in $x,y${\rm)} and given by
\begin{equation}
Q(x,y;z)=\sum_{i,j,k=0}^{\infty}q_{i,j,k}x^{i}y^{j}z^{k}=\frac{c_0}{(1-\alpha_0 x)(1-\beta_0 y)}.
\end{equation}
\end{prop}

\begin{proof}
Let us first construct all pairs $(\alpha,\beta)$ such that $c\alpha^i\beta^j$ solves \eqref{eq_i>0_j>0_sumS}, \eqref{eq_i>0_j=0_sumS} and \eqref{eq_i=0_j>0_sumS}. We obtain that $\alpha$ and $\beta$ must satisfy
\begin{align}
\alpha\beta/z&=\beta+\beta^2+\alpha\beta^2+\alpha^2\beta+\alpha^2\beta^2, \label{try1}\\
\alpha\beta/z&=\alpha+\alpha\beta^2+\alpha^2\beta+\alpha^2\beta^2. \label{try2}
\end{align}
(Note that the reason why we obtain two and not three equations is that Equation \eqref{eq_i>0_j=0_sumS} gives exactly the same information as \eqref{eq_i>0_j>0_sumS}.) Further, all pairs $(\alpha,\beta)$ that satisfy both \eqref{try1} and \eqref{try2} are such that $\beta(1+\beta)=\alpha$. Substituting the latter into \eqref{try1} gives \eqref{rel:bet}. Denote by $\beta_k$, $k\in\{0,\ldots ,3\}$, the four roots of \eqref{rel:bet}. Among the $\beta_k$, only one, say $\beta_0$, is a unique formal power series in $z$, see \cite[Proposition 6.1.8]{STAN}. Letting
     \begin{equation}
          \label{albe}
          \alpha_0=\beta_0(1+\beta_0),
     \end{equation}
we conclude that for any $c_0$, $c_0\alpha_0^i\beta_0^j$ is a solution to \eqref{eq_i>0_j>0_sumS}, \eqref{eq_i>0_j=0_sumS} and \eqref{eq_i=0_j>0_sumS}. With the particular choice $c_0=1/[1-z(\alpha_0+\alpha_0\beta_0+\beta_0)]$, $c_0 \alpha_0^i\beta_0^j$ is also solution to  \eqref{eq_i=0_j=0_sumS}. By uniqueness of the solution to the recursion relations (see Lemma \ref{lemm1}), we conclude that $q_{i,j}=c_0\alpha_0^i\beta_0^j$. Finally, by \eqref{rel:bet} and \eqref{albe}, we can simplify $c_0$ into $\beta_0/z$.
\end{proof}

We now modify the boundary behavior of this walk and we show that this has severe consequences for the GF.  Consider again the walk
with $\mathcal{S}=\{(-1,0),(-1,-1),(0,-1),$
$(1,-1),(1,0)\}$ in the interior, but this time with
$\mathcal{S}_V=\{(0,-1),(1,-1),(1,0)\}$ on the vertical boundary, $\mathcal{S}_0=\{(1,0)\}$, and the following rather special step set on the horizontal boundary: if the walk is in state $(i,0)$ with $i>0$, the possible steps are $\{(-1,0),(1,0)\}$ and a {\it big} step to $(0,i)$;
see the right display in Figure \ref{other_examples}. We should note that a {\it random} walk, with similar unusual boundary behavior, has already been considered in \cite{MR2042259}.
The equations for $q_{i,j}$ are then given by \eqref{eq_i>0_j>0_sumS}, \eqref{eq_i>0_j=0_sumS}, \eqref{eq_i=0_j=0_sumS} and
\begin{align}
          \label{eq_i=0_j>0_sumS2}
          q_{0,j}/z&=q_{j,0}+q_{0,j+1}+q_{1,j}+q_{1,j+1},
          \quad j>0.
\end{align}

\begin{prop}
\label{partcase}
Let
\begin{equation}
\label{innerss}
g(t)=\frac{1-tz-\sqrt{(1-tz)^2-4z^2(1+t)^2}}{2z(1+t)},
\end{equation} $\beta_0=0$,
and define $c_0$ by
     \begin{equation}
     \label{cft}
          1/c_0={\displaystyle\sum_{k=0}^{\infty} \frac{\beta_2\cdots \beta_{k+1}}{(1+\beta_1)\cdots (1+\beta_k)}[1-z(\beta_k+\beta_{k+1}+\beta_k\beta_{k+1})]}
     \end{equation}
with
     \begin{equation}
     \label{innersss}
          \beta_{k+1}=\alpha_{k},
          \quad \alpha_{k+1}=g(\beta_{k+1}),
          \quad c_{k+1}=c_{k}\frac{\alpha_{k+1}}{1+\beta_{k+1}},\quad  k\geq 0.
     \end{equation}
The unique solution to Equations \eqref{eq_i>0_j>0_sumS}, \eqref{eq_i>0_j=0_sumS}, \eqref{eq_i=0_j=0_sumS} and  \eqref{eq_i=0_j>0_sumS2} is given by
     \begin{equation}
     \label{sdsd1s}
          q_{i,j}=\sum_{k=0}^\infty c_k \alpha_k^i\beta_k^j.
     \end{equation}
Therefore, the GF is  given by
     \begin{equation}
          Q(x,y;z)=\sum_{k=0}^{\infty} \frac{c_k}{(1-\alpha_k x)(1-\beta_k y)},
     \end{equation}
which has infinitely many poles in $x$ and $y$, and is therefore nonholonomic as a trivariate function.
\end{prop}

\begin{proof}
Let us first notice that compared to Section \ref{The_compensation_approach}, there is here a main novelty in the application of the compensation approach: it comes from the fact that Equation \eqref{eq_i>0_j>0_sumS} for the interior and Equation \eqref{eq_i>0_j=0_sumS} for the horizontal boundary are exactly similar. Accordingly, the compensation terms will appear by repeatedly considering Equation \eqref{eq_i=0_j>0_sumS} for the vertical boundary only.

We start from the product $q_{i,j}=c\alpha^i\beta^j$. The pairs $(\alpha,\beta)$ that satisfy the inner equations \eqref{eq_i>0_j>0_sumS} are characterized by $\alpha=g(\beta)$ with $g$ as in \eqref{innerss}. (It is worth mentioning that a priori, we have two possibilities for $g$: one as in \eqref{innerss}, and another with a sign $+$ in front of the square root. We have to choose \eqref{innerss} because otherwise, the $\alpha_k$ and $\beta_k$ in \eqref{innersss} would not be power series.) Substituting $c\alpha^i\beta^j$ into
\eqref{eq_i=0_j>0_sumS2} yields
\begin{equation}\label{r15}
c\alpha^j=c\beta^j(1/z-\alpha-\beta-\alpha\beta).
\end{equation}
This suggests to set $\beta_0 = 0$, as otherwise \eqref{r15} does not hold. (It is worth noting that this choice of $\beta_0$ is less canonical than the choice we made in Lemma \ref{lem_a_b}. However, if we are able to construct a solution to the recursion relations \eqref{eq_i>0_j>0_sumS}, \eqref{eq_i>0_j=0_sumS}, \eqref{eq_i=0_j=0_sumS} and  \eqref{eq_i=0_j>0_sumS2} with this particular choice, then by uniqueness of the solution, it will be the only one, see Lemma \ref{lemm1}.)
Consider thus the product $c_0\alpha^i_0\beta^j_0$, with $\beta_0=0$ and $\alpha_0=g(\beta_0)$, for which \eqref{r15} becomes
\begin{equation}\label{r16}
c_0\alpha_0^j=0.
\end{equation}
We then conclude that $c_0\alpha^i_0\beta^j_0$ does not satisfy \eqref{eq_i=0_j>0_sumS2} and to correct for the error, we add a second product $c_1\alpha^i_1\beta^j_1$, with $(\alpha_1,\beta_1)$ on the curve $\alpha=g(\beta)$, in such a way that the term $c_1\alpha^i_1\beta^j_1$ appears on the right-hand side of \eqref{r16}. This yields
\begin{equation}\label{r17}
c_0\alpha_0^j=c_1\beta_1^j(1/z-\alpha-\beta-\alpha\beta),
\end{equation}
which gives
\begin{equation}
\beta_1=\alpha_0,\quad \alpha_1=g(\beta_1),\quad c_1=c_0\frac{1}{1/z-\alpha_1-\beta_1-\alpha_1\beta_1}=c_{0}\frac{\alpha_1}{1+\beta_1},
\end{equation}
where the last equality follows from \eqref{eq_i>0_j>0_sumS}. By adding the second product, we also obtain an extra term $c_1\alpha_1^j$ on the left-hand side of \eqref{r17}. That is why we add a third product $c_2\alpha^i_2\beta^j_2$, with $(\alpha_2,\beta_2)$ on the curve $\alpha=g(\beta)$, and via similar reasoning,
\begin{equation}
\beta_2=\alpha_1,\quad \alpha_2=g(\beta_2),\quad c_2=c_{1}\frac{\alpha_2}{1+\beta_2}.
\end{equation}
Continuing this procedure means that we keep adding products $c_k\alpha^i_k\beta^j_k$, with $c_k$, $\alpha_k$ and $\beta_k$ as in \eqref{innersss}, so that (\ref{sdsd1s}) is a formal solution to Equations \eqref{eq_i>0_j>0_sumS}, \eqref{eq_i=0_j=0_sumS} and  \eqref{eq_i=0_j>0_sumS2}.

We now find $c_0$. With \eqref{innersss}, that we have proved above, we obtain
\begin{equation}
     \frac{c_k}{c_0} = \frac{\beta_2\cdots \beta_{k+1}}{(1+\beta_1)\cdots (1+\beta_k)},
     \quad k\geq 0.
\end{equation}
Further, with \eqref{eq_i=0_j=0_sumS} we have
\begin{equation}
     \frac{c_0}{z}\sum_{k=0}^{\infty}\frac{c_k}{c_0}=\frac{1}{z}+c_0\sum_{k=0}^{\infty}\frac{c_k}{c_0}(\alpha_k+\beta_k+\alpha_k\beta_k).
\end{equation}
Together with \eqref{innersss}, we obtain \eqref{cft}.
The $\alpha_k$ and $\beta_k$ converge to the unique formal power series $\beta_*$ solution to $g(\beta_*)=\beta_*$. Using \eqref{innerss}, we obtain that
     \begin{equation}
     \label{minimal_poly_beta}
          \beta_*^3+2\beta_*^2+\beta_*(1-1/z)+1=0.
     \end{equation}
The series $\beta_*$ is such that $\beta_*(0)=0$. It remains to show that $q_{i,j}$ defined in \eqref{sdsd1s} is a (well-defined) formal power series. There is here an interesting phenomenon: for the general compensation approach of Section \ref{The_compensation_approach}, $q_{i,j}$ was a formal power series because $\alpha_k=z^{2k+1}\widehat{\alpha}_k$  and
$\beta_k=z^{2k+2}\widehat{\beta}_k$, where $\widehat{\alpha}_k$ and $\widehat{\beta}_k$
are formal power series (see Proposition \ref{prop_control}). But here, we have $\alpha_k =z\widehat{\alpha}_k$ and
$\beta_k=z\widehat{\beta}_k$, where  $\widehat{\alpha}_k$ and $\widehat{\beta}_k$ are formal power series not vanishing at $0$ (this follows from \eqref{innerss} and \eqref{innersss}). On the other hand, we have $c_k = z^k \widehat{c}_k$, where $\widehat{c}_k$ is a formal power series. Therefore, in the particular case of Proposition \ref{partcase}, $q_{i,j}$ is a formal power series thanks to the constants $c_k$ (see \eqref{innersss} and \eqref{sdsd1s}).
%
%
The nonholonomy of the trivariate function $Q(x,y;z)$ follows by the same arguments as that given below Proposition \ref{infinite_poles}.
\end{proof}

\noindent{\textbf{Acknowledgments.}}
The work of Johan S.H.\ van Leeuwaarden was supported by an ERC Starting Grant.
Kilian Raschel would like to thank EURANDOM for providing wonderful working conditions during two visits in the year $2010$. His work was partially supported by CRC $701$, Spectral Structures and Topological Methods in Mathematics at the University of Bielefeld. The authors thank an anonymous referee for numerous valuable and constructive comments.

\appendix{}

\section{Remaining proofs}
\label{fa}

Let us start by expressing $h$, defined in (\ref{def_h}), as an alternating sum.
The reason why we wish to formulate $h$ differently is twofold. First, this will enable us to
show that $\rho\in[0.34499975,0.34499976]$, see Corollary \ref{approx_rho}---and
in fact, in a similar way, we can approximate $\rho$ up to any level of precision. Also,
this is actually a key lemma for proving Lemma \ref{def_rho}.
\begin{lem}
\label{other_form}
     \begin{equation}
     \label{h_T}
          h(z)=1+2z\bigg(-1-\alpha_0+\sum_{k=0}^{\infty}(-1)^kT_k\bigg),
          \quad T_k=\left\{\begin{array}{lcl}
          \beta_{k/2}\alpha_{k/2} &  \ {\rm if} \  k  \ {\rm is  \ even},\\
          \beta_{(k-1)/2}\alpha_{(k+1)/2} & \ {\rm if} \  k \ {\rm is \  odd}.
          \end{array}\right.
     \end{equation}
\end{lem}

\begin{proof}
Equation \eqref{h_T} follows immediately from \eqref{nfx} and \eqref{def_h}.
\end{proof}

As a preliminary result, we also need the following refinement of Proposition \ref{prop:analytic-properties}:

\begin{lem}
\label{extension}
Assume that $z\in(0,1/\sqrt{8})$. Then both sequences $\{\alpha_k\}_{k\geq 0}$ and
$\{\beta_k\}_{k\geq 0}$ are positive and decreasing. Moreover, for all $k \geq 0$
     \begin{equation}
          0\leq \alpha_k\leq 1/\sqrt{2}^{2k+1},
          \quad 0\leq \beta_k\leq 1/\sqrt{2}^{2k+2}.
     \end{equation}
\end{lem}
\begin{proof}
It is enough to prove that for all $\vert t\vert \leq 1$ and $z\in(0,1/\sqrt8)$,
\begin{equation}
     0\leq f(t)\leq \frac{t}{\sqrt2}.
\end{equation}
The proof of this inequality is similar to  that of \eqref{properties_f} and therefore omitted.
\end{proof}

\begin{cor}
\label{approx_rho}
The series $h(z)$ defined in \eqref{def_h} or \eqref{h_T} has a radius of convergence at least $1/\sqrt{8}$. Its first positive zero $\rho$ lies in $[0.34499975,0.34499976]$.
\end{cor}
\begin{proof}
For $z\in (0,1/\sqrt{8})$, the sequence $\{T_k\}_{k\geq 0}$ defined in Lemma \ref{other_form}
is positive and decreasing, see Lemma \ref{extension}. In
particular, denoting
     \begin{equation}
     \label{def_approx}
          \Lambda^{p}=\sum_{k=0}^{p}(-1)^kT_k
     \end{equation}
and using \eqref{h_T}, for all $p\geq 0$ and all $z\in (0,1/\sqrt{8})$ we have
     \begin{equation}
          1+2z(-1-\alpha_0+\Lambda^{2p+1})<
          h(z)<1+2z(-1-\alpha_0+\Lambda^{2p}).
     \end{equation}
Applying the last inequalities to $p=4$ and noting that the right-hand side
(resp.\ left-hand side) evaluated at $0.34499976$ (resp.\ $0.34499975$) is negative
(resp.\ positive) concludes the proof.
\end{proof}

\begin{proof}[Proof of Lemma \ref{def_rho}]
First, using the explicit expression of $\alpha_k$ and $\beta_k$, and employing calculus software, we obtain that the algebraic function $1+2z(-1-\alpha_0+\Lambda^{5})$ has only one zero within the circle of radius $1/\sqrt{8}$, and that
                \begin{equation}
                     \label{min_4}
                     \inf_{|z|= 1/\sqrt{8}}|1+2z(-1-\alpha_0+\Lambda^{5})|>10^{-2}.
                \end{equation}
Note that the lower bound $10^{-2}$ in \eqref{min_4} is almost optimal; it is rather small because like $\rho$, the unique zero of $1+2z(-1-\alpha_0+\Lambda^{5})$ within the disc of radius $1/\sqrt{8}$ is very close to $1/\sqrt{8}$.

By Rouch\'e's theorem (see Theorem \ref{RTT}), applied to the circle with radius $1/\sqrt{8}$, it is now enough to prove that
                \begin{equation}
                     \label{min_infinity}
                     \sup_{|z| = 1/\sqrt{8}}\left|[1+2z(-1-\alpha_0+\Lambda^{5})]-[1+2z(-1-\alpha_0+\Lambda^{\infty})]\right|<10^{-2}.
                \end{equation}
For this we write
     \begin{equation}
     \label{mup}
          \left|\Lambda^{\infty}-\Lambda^p\right|=\bigg|\sum_{k=p+1}^{\infty}(-1)^kT_k\bigg|
          \leq \sum_{k=p+1}^{\infty}|T_k|\leq \frac{1/\sqrt{2}^{2p+5}}{1-(1/\sqrt{2})^{2}}< 1/\sqrt{2}^{2p+3},
     \end{equation}
where the second upper bound is obtained from the inequality $|T_k|\leq 1/\sqrt{2}^{2k+3}$, see
\eqref{h_T} and Lemma \ref{extension}. By \eqref{mup} we then obtain
     \begin{equation}
          \left|[1+2z(-1-\alpha_0+\Lambda^{p})]-[1+2z(-1-\alpha_0+\Lambda^{\infty})]\right|\leq 2|z|\frac{1}{2^{p+3/2}}.
     \end{equation}
If $|z|< 1/\sqrt{8}$, the last quantity is bounded from above by $1/2^{p+2}$. For $p=5$, we obviously get $1/2^{p+2}<10^{-2}$,
and \eqref{min_infinity} is proven.
\end{proof}

\begin{lem}
\label{approx_mul}
     \begin{equation}
          C_{0,0}=0.0531\cdot [1+O(10^{-3})].
     \end{equation}
\end{lem}

\begin{proof}
As in the proof of Proposition \ref{prop_asymp_ij}, we can write
     \begin{equation}
          C_{0,0}=\frac{3\rho -1}{\rho[1-\rho^{2}\widehat{x}_{0,0}'(\rho)]},
     \end{equation}
so that the main difficulty lies in approximating
$\widehat{x}_{0,0}'(\rho)$. For this, we use \eqref{nfx} and \eqref{h_T} to write
     \begin{equation}
     \label{eq_p}
          \widehat{x}_{0,0}'(z)=-2\alpha_{0}'(z)+\sum_{k=0}^{\infty}(-1)^{k}T_{k}'(z).
     \end{equation}
While it is easy to control the series $\sum_{k=0}^{\infty}(-1)^{k}T_{k}(z)$, because
$\alpha_k(z)$ and $\beta_k(z)$, and hence $T_{k}(z)$,
decrease exponentially fast to $0$ as $k\to \infty$, see Proposition
\ref{prop_control} and Lemma \ref{extension}, it is not obvious how we should deal with $\sum_{k=0}^{\infty}(-1)^{k}T_{k}'(z)$,
where terms like $\alpha_{k}'(z)$ and $\beta_k'(z)$ appear. We next show
that $\alpha_{k}'(z)$ and $\beta_k'(z)$ actually also decrease exponentially fast to $0$ as $k\to \infty$,
at least for $z\in[1/4,0.35]$. Note that this assumption on $z$ is not restrictive, since $\rho$
belongs to the interval $[1/4,0.35]$ by Corollary \ref{approx_rho}.

Consider the sequence $\{\gamma_k(z)\}_{k\geq 0}$ defined by $\gamma_0(z)=\alpha_0(z)$ and,
for $k\geq 0$, by $\gamma_{k+1}(z)=f(\gamma_k(z))$, with $\alpha_0(z)$ and $f$ as defined in \eqref{def_alpha_0}
and \eqref{properties_alpha_beta}, respectively. Note that \eqref{properties_alpha_beta} gives
$\gamma_{2k}(z)=\alpha_k(z)$ and $\gamma_{2k+1}(z)=\beta_k(z)$. The sequence $\{\gamma_k'(z)\}_{k\geq 0}$
satisfies the recurrence relation
     \begin{equation}
     \label{rr}
          \gamma_{k+1}'(z)=\gamma_k'(z)\partial_t f(\gamma_k(z))+\partial_z f(\gamma_k(z)).
     \end{equation}
Below we study the coefficients $\partial_t f(\gamma_k(z))$ and $\partial_z f(\gamma_k(z))$
of \eqref{rr}.
By using expression \eqref{properties_alpha_beta} of $f$, we easily obtain that
for all $z\in[0,1/\sqrt{8}]$ and all $t\in[0,a_{0}(z)]$
     \begin{equation}
     \label{m1}
          |\partial_t f(\gamma_k(z))|\leq 4\sqrt{2}/9.
     \end{equation}
In addition, $\partial_z f(t)=f(t)/[z\sqrt{1-4z^2(1+t^{2})}]$, in such a way that
$|\partial_z f(t)|\leq |2/z|\cdot |f(t)|$ for all $z\in[0,1/\sqrt{8}]$ and all
$t\in[0,a_{0}(z)]$. Using then Lemma \ref{extension},
we obtain
     \begin{equation}
     \label{m2}
          |\partial_z f(\gamma_{k}(z))|\leq |2/z|\cdot |f(\gamma_{k}(z))|=|2/z|\cdot |\gamma_{k+1}(z)|
          \leq |2/z| \cdot 1/\sqrt{2}^{k+1}\leq 8/\sqrt{2}^{k+1},
     \end{equation}
where the last inequality follows from the assumption $z\in[1/4,0.35]$.
From \eqref{rr}--\eqref{m2} we get
     \begin{equation}
          |\gamma_{k+1}'(z)|\leq \eta \cdot |\gamma_{k+1}'(z)|+\xi_k,
          \quad \eta=4\sqrt{2}/9,
          \quad \xi_k=8/\sqrt{2}^{k+1}.
     \end{equation}
We deduce that
     \begin{equation}
     \label{mmaa}
          |\gamma_{k+1}'(z)|\leq \eta^{k+1}\cdot |\gamma_{0}'(z)|+
          \sum_{p=0}^{k}\eta^p \xi_{k-p}\leq \eta^{k+1} \cdot |\gamma_{0}'(z)|+{72}/{\sqrt{2}^{k+1}}
          \leq 100/\sqrt{2}^{k+1},
     \end{equation}
where the last inequality follows from $\sup_{z\in[0.25,0.35]}|\gamma_{0}'(z)|\leq 13$
and $\eta\leq 1/\sqrt{2}$. With Lemma \ref{extension} and \eqref{mmaa},
we obtain
\begin{equation}
|T_k'(z)|\leq 200/\sqrt{2}^{k+2},
\end{equation} and thanks to \eqref{eq_p} we can obtain approximations of $\widehat{x}_{0,0}'(\rho)$
up to any level of precision. Lemma~\ref{approx_mul} follows.
\end{proof}
\end{document}